# Calculating the natural density of Mersenne numbers using nonstandard mathematical analysis.

Gennady Eremin

**Abstract.** Currently, among the open (unsolved) problems in number theory is the following: it is unknown what is the natural density of the sequence of Mersenne numbers in the set of natural numbers. In the paper, using methods of nonstandard mathematical analysis, we obtain the following equation: the natural density of Mersenne numbers (some infinitesimal value $\varepsilon$) multiplied by the sum of the reciprocals of odd numbers (the infinitely large value $\omega = 1 + 1/3 + 1/5 + ...$) is equal to 1, or the equality $\varepsilon = 1/\omega$ is true. In nonstandard analysis, the resulting infinitesimal numbers $\varepsilon$ and $1/\omega$ are considered equivalent. We obtained this result by working with a two-dimensional matrix of non-negative integers, where odd numbers are separated from even ones by the Pepis-Kalmar pairing function.

# 1 Introduction

There is an open (unsolved) problem in number theory, namely: what is the natural density of the sequence of Mersenne numbers in the set of natural numbers, let's clarify, in the set of non-negative integers $\mathbb{N} = \{0, 1, 2, 3, ...\}$. Is the density of Mersenne numbers zero or positive (greater than zero)? In this paper, we solve this problem by working with a two-dimensional matrix of natural numbers; methods and techniques of nonstandard mathematical analysis are used to calculate the density of Mersenne numbers.

**1.1. Matrix of natural numbers.** In the matrix of natural numbers $\mathbb{M}$, odd numbers are separated from even ones using the Pepis-Kalmar pairing function [2]. The figure shows the initial terms of the matrix, on the $x$-axis we have even numbers starting from 0, and on the $y$-axis we have Mersenne numbers (see the sequence A000225 in [6]). In the columns of the matrix, each current odd term is generated by the previous term of the same column according to the recurrence formula of Mersenne numbers: a term $t$ generates the subsequent (and only one) term $t'$ of the same column, i.e. $t' = 2t + 1$. It is logical to think of the columns as unary trees, or *Mersenne trees*.

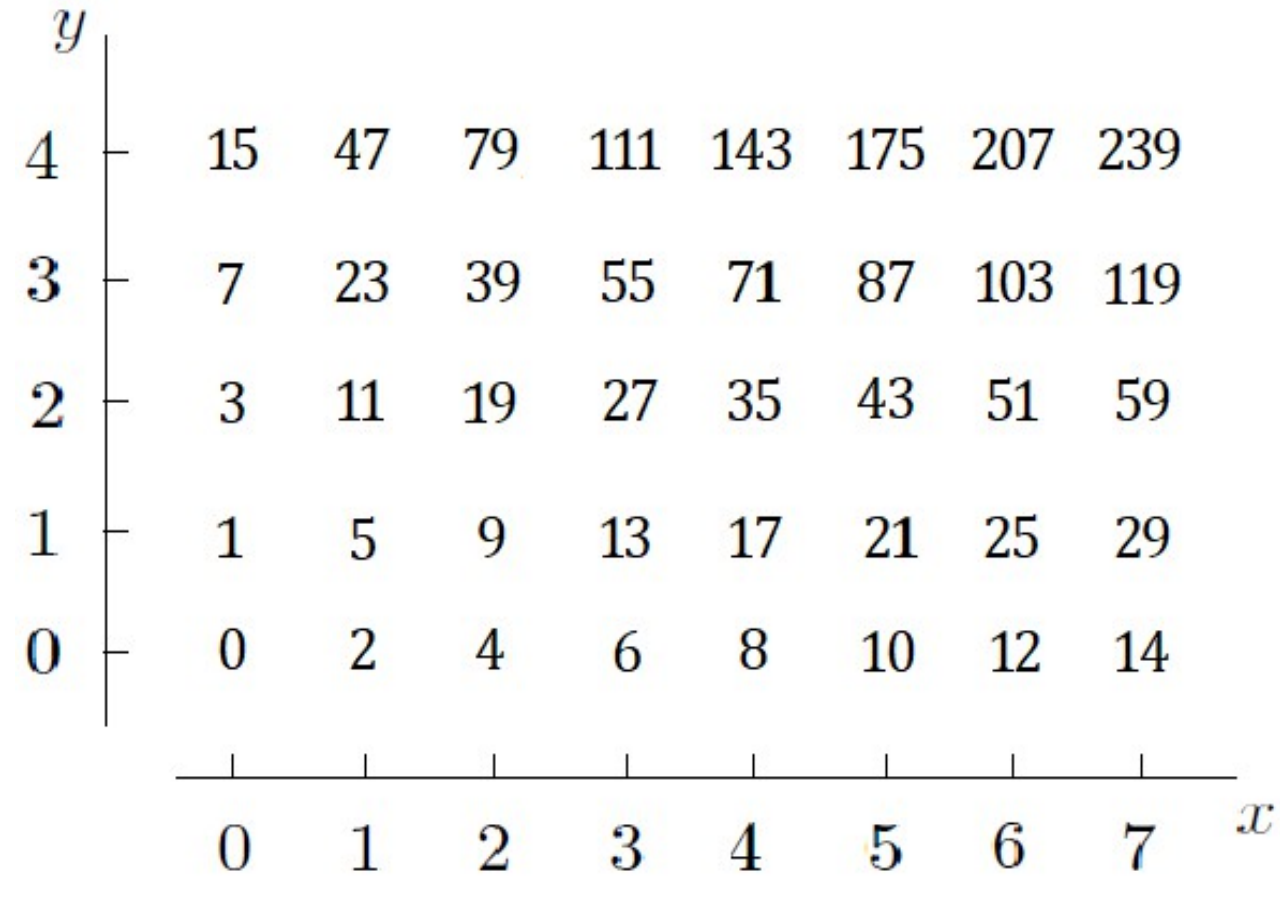

Each column of the matrix $\mathbb{M}$ starts with an even number (*tree root*). The set of even natural numbers is infinite, so we have an infinite number of unary Mersenne trees. Mersenne trees are disjoint because not only the successor formula is linear, but also the precedence formula $t = (t' - 1)/2$. In Mersenne trees, only roots are even numbers.

The initial column of the matrix is the Mersenne numbers 0, 1, 3, 7, 15, ... The binary representation of the $n$-th Mersenne number $2^n - 1$ is $n$ units. The next Mersenne tree is the Thabit numbers 2, 5, 11, 23, 47, ... (see the sequence A055010 in [6]). The binary representation of the $n$-th Thabit number $3 \times 2^n - 1$ is $n+2$ digits long, consisting of "10" followed by $n$ 1s ($n$ binary units). Here is the general formula for the matrix term (the Pepis-Kalmar pairing function):

(1) $$f(x, y) = (2x+1) \times 2^y - 1.$$

Obviously, in formula (1), the power of two $y$ is the number of final units in the binary code of the resulting odd number (lines of the matrix with coordinate $y > 0$). In case $y = 0$ we get even numbers on the $x$-axis. From formula (1) it is easy to obtain a recurrence formula in the columns of the matrix: for an arbitrary term $t = f(x, y)$, the subsequent term in the same column is equal to

$$t' = f(x, y+1) = (2x+1) \times 2^{y+1} - 1 = 2(2x+1) \times 2^y - 2 + 1 = 2t + 1.$$

Let us repeat, the matrix $\mathbb{M}$ contains all non-negative integers, i.e. all elements of the set $\mathbb{N}$. Let's prove the following statement.

**Proposition 1.** *In the matrix of natural numbers $\mathbb{M}$, terms do not repeat.*

*Proof.* Let's start with the $x$-axis, which contains all even natural numbers in ascending order, starting with zero. Even numbers do not repeat, since on the $x$-axis we have an arithmetic progression with a step of 2. Accordingly, the binary codes of even numbers do not repeat either. And then in the columns, with each step (stepwise increase of the $y$-coordinate by one), in each term we add unit at the end of the binary code, which is equivalent to first multiplying by 2 (0 is added at the end of the binary code), and then adding 1 (while at the end of the binary code zero changes to unit). As a result, the uniqueness of binary codes is respected, and therefore there are no repeating terms in the column family of the matrix, which contains all natural numbers (including 0). □

**1.2. A bit on nonstandard mathematical analysis.** To calculate the natural density of Mersenne numbers, we use nonstandard analysis, which differs from the usual standard analysis by adding infinitesimal (non-zero) and infinitely large quantities [1, 3, 4]. The following is important: in nonstandard analysis, infinitesimals are considered not as variable quantities (or as functions tending to zero), but as constant actual numbers (albeit imaginary, ideal ones). Every non-zero infinitesimal quantity is so close to 0 that when we multiply it by any finite (standard) number, we get an infinitesimal quantity as a result again. Additionally, we note that 0 is also considered an infinitesimal number, and this is logical.

With infinitesimal quantities and standard numbers, you can perform ordinary arithmetic operations (addition, multiplication, etc.). Back in the 17th century, Gottfried Leibniz (a German mathematician and philosopher) formulated the rules for operating with infinitesimals in the form of calculus. For example, when multiplying a non-zero infinitesimal by another non-zero infinitesimal (or squaring an infinitesimal), we obtain a new non-zero infinitesimal of higher order, which is even closer to zero. We obtain an infinitely large number if we divide a finite standard number by a non-zero infinitesimal (dividing by zero is prohibited in both standard and nonstandard analysis). The converse is also true: if we divide a non-zero standard number by an infinitely large number, we get a non-zero infinitesimal quantity.

In nonstandard analysis, in every infinite set of objects there is at least one nonstandard element (*idealization principle*). For example, in our matrix, each line and each column end with a unique infinitely large value. In this paper we obtain an infinite number of unique infinitesimal quantities whose sum is equal to 1.

# 2 Natural density in matrix lines.

In number theory, the natural density of an arbitrary set of natural numbers characterizes the size of such a set, i.e. how large the given set is relative to the set of all natural numbers $\mathbb{N}$ (in this case, the natural density of the entire set $\mathbb{N}$ is taken to be equal to 1). For example, the set of even natural numbers has a natural density of 1/2, since every second natural number is even; the set of odd natural numbers also has the same density. The density of every finite set of natural numbers is 0 (and this is logical, since the set $\mathbb{N}$ is infinite). Obviously, the natural density of the set $\mathbb{N}$ and, accordingly, the density of the terms of the matrix $\mathbb{M}$ are the same and equal to 1. Accordingly, the total density of all matrix lines is also equal to 1.

In each line of the matrix $\mathbb{M}$ we have an infinite arithmetic progression, and there are an infinite number of such progressions. If we fix the $y$-coordinate in formula (1), we obtain an arithmetic progression of the form $a + dx$, where the initial term $a = 2^y - 1$, and the difference in the progression $d = 2^{y+1}$. In every arithmetic progression, the initial term is the corresponding Mersenne number, the difference of the progression is a fixed power of two, starting with 2, 4, 8, ... As is known, the natural density of such progressions is determined quite simply and is equal (starting from the bottom line, $x$-axis) to 1/2, 1/4, 1/8, ..., respectively, i.e., the values inverse to the differences of the progressions. The sum of such an infinitely decreasing geometric progression is well known

$$S = 1/2 + 1/4 + 1/8 + 1/16 + ... = 1,$$

as a result we obtain the natural density of the set $\mathbb{N}$. So, with the lines of the matrix, everything is quite simple and obvious; next, we will work with the columns of the matrix, and unfortunately, everything is more complex and not so straightforward in the columns.

# 3 Natural density in matrix columns.

In the matrix of natural numbers, terms are not repeated, therefore the density of all columns, as well as the density of all lines, are the same and equal to 1. The initial column contains Mersenne numbers (terms with zero $x$-coordinate, i.e., $y$-axis). Since the $n$-th Mersenne number is $2^n - 1, n \geq 0$, we can state the following: in the range of natural numbers $\{0, 1, 2, 3, ..., 2^n - 1\}$ (the length of such a range is $2^n$) we have Mersenne numbers with indices $0, 1, 2, ..., n$ (i.e., a total of $n + 1$ Mersenne numbers). Thus, in the set of natural numbers $\mathbb{N}$ in the initial range of length $2^n$, the density of Mersenne numbers is $(n+1)/2^n$. As $n$ increases, the density of Mersenne numbers tends to zero (in the resulting fraction it is enough to compare the first derivatives of the numerator and the denominator, we use the well-known L'Hopital rule). As a result, we obtain an infinitesimal quantity, let's denote this number as $\varepsilon$. Let us remember that in nonstandard analysis, infinitesimals are not variables or functions tending to zero, but special actual numbers, often different from zero. In this case, it is theoretically possible that $\varepsilon = 0$, and we need to prove the positive natural density of Mersenne numbers, i.e. $\varepsilon > 0$.

The next matrix column (coordinate $x = 1$) contains Thabit numbers $f(1, y) = 3\times 2^y - 1$. The natural density of such numbers is three times less than the density of Mersenne numbers, since $n+1$ Thabit numbers are located in the range $\{0, 1, 2, 3, ..., 3\times 2^n - 1\}$, the length of such a range is $3\times 2^n$ (three times greater than the corresponding range with the same number of Mersenne numbers). Since in nonstandard analysis one can perform arithmetic operations with infinitesimals, it is logical to accept that the natural density of Thabit numbers is equal to $\varepsilon/3$ (we reduce the density of Mersenne numbers by a factor of three). Let us add that in nonstandard analysis such two infinitesimal quantities $\varepsilon$ and $\varepsilon/3$ are called infinitesimals of the same order.

Next, for the matrix column with terms of the form $f(2, y) = 5\times 2^y - 1$ (coordinate $x = 2$) we obtain another infinitesimal density equal to $\varepsilon/5$. Thus, three unique infinitesimals of the same order were formed: $\varepsilon$, $\varepsilon/3$ and $\varepsilon/5$. Further, in the following columns we obtain infinitesimals of the same order $\varepsilon/7$, $\varepsilon/9$ and so on. As a result, we have a sequence of infinitesimal quantities of the same order, which gradually approach 0, and there are obviously a countable number of such infinitesimal numbers. In the resulting sequence, all infinitesimals are non-zero and unique (but this is in the case $\varepsilon > 0$).

Thus, near zero we obtain a family of positive infinitesimal quantities, each of which is infinitely close to zero. In nonstandard analysis, such a family of infinitesimals is some part of the *monad* of the number 0. (Let us explain: a monad is a vicinity of a point on the number line; the monad includes the point itself and all infinitesimal deviations from it.)

# 4 Natural density of Mersenne numbers.

In the columns of the matrix we obtain sequences of natural numbers with infinitely small natural density. Each such infinitesimal quantity is unique. Since in the matrix $\mathbb{M}$ the total natural density of all columns (as well as the total density of lines) is equal to the density of the set of natural numbers $\mathbb{N}$, we obtain the following equality:

(2) $$1 = \varepsilon + \varepsilon/3 + \varepsilon/5 + \varepsilon/7 + ... = \varepsilon \times (1 + 1/3 + 1/5 + 1/7 + ...) = \varepsilon \times \omega.$$

From equality (2) it follows that $\varepsilon > 0$ (other terms of the sum are also positive), otherwise each term and the sum itself will be equal to zero. As a result, we obtain the following: one is equal to the product of the Mersenne number density $\varepsilon$ and an infinitely large value $\omega$, which is the sum of all reciprocals of odd natural numbers. Let us formulate the corresponding statement.

**Theorem 2.** *The natural density of Mersenne numbers and the reciprocal of the sum of reciprocals of odd natural numbers are equivalent infinitesimal numbers, i.e.* $\varepsilon = 1/\omega$ or $\varepsilon \times \omega = 1$.

Let us repeat that in mathematics the distribution density of Mersenne numbers in the set of natural numbers is still unknown. In the paper [2] we proved that the density of Mersenne numbers is positive (non-zero). Now we have obtained a concrete equation for determining the natural density of a sequence of Mersenne numbers.

**Acknowledgements.** The author is very grateful to Olena Kachko for discussion and useful comments on preliminary versions of the paper. Also the author would like to thank the anonymous referees for their valuable comments and suggestions.

*Email address*: ergenns@gmail.com

Written: July 7, 2026.